\documentclass[11pt]{article}
\usepackage{amssymb}
\usepackage{graphicx}
\usepackage{amsmath}
\usepackage{multirow}
\newtheorem{theorem}{Theorem}

\newtheorem{lemma}[theorem]{Lemma}

\begin{document}
\title{{\bf On some combinatorial formulae coming from Hessian Topology}}
\author{ 
Adriana Ortiz-Rodr\'{\i}guez and Federico S\'anchez-Bringas}
\maketitle

\begin{abstract}
The interaction between combinatorics and algebraic and differential geometry is
very strong. While researching a problem of Hessian topology, 
we came across a series of identities 
of binomial coefficients, which are
useful for proving a topological property of 
certain spaces whose elements are graphs of a class
of hyperbolic polynomials.  
These identities are proven by different methods
in combinatorics.
\end{abstract}


\noindent {\small {\it MS classification. 05A19, 57R52, 53A05}}

\section{Introduction}\label{intro}
Examples of the interaction between combinatorics and algebraic geometry are
abundant:  the polynomial method that was recently surveyed by
Terence Tao \cite{tao}; the Combinatorial Nullstellensatz of Noga Alon \cite{alon};
the affirmative answer to the conjecture of Read and Rota-Heron-Welsh by
June Huh \cite{huh}; and, Matthew Baker and Serguei Norine's  graph-theoretic
analogue of the classical Riemann-Roch theorem \cite{baker}.

All of these examples seem to indicate that the interaction is one way. 
This, of course,  is not true. However, the use of combinatorics in the 
proof of theorems in algebraic geometry is more subtle. Examples of combinatorial 
arguments in algebraic geometry are: computing maximum possible number of ordinary 
double points for a surface of degree 5 by Arnaud Beauville \cite{beauville}; 
giving an explicit combinatorial formula for the structure constants
of the Grothendieck ring of a Grassmann variety with respect to its basis
of Schubert structure sheaves by Anders Skovsted Buch \cite{buch}; and, the 
study of Donaldson-Thomas invariants by Benjamin Young \cite{young}.
The interaction of differential topology and geometry with combinatorics is even 
more obvious due to discrete Morse theory~\cite{forman_98} and combinatorial differential geometry 
by Robin Forman \cite{forman_99}; and, lattice gauge theory by Kenneth G. Wilson~\cite{wilson}.

In this article, we prove as a result of some combinatorial identities that 
the graphs of the class of hyperbolic polynomials studied in  \cite{geoast}
have a topological property which guarantees the existence of a number, depending on the degree, 
of different connected components of hyperbolic homogeneous polynomials. This fact is relevant for
the description of the Hessian topology of hyperbolic polynomials. 

The paper is organised as follows: In Section~\ref{prelim}, we give a detail 
account of the problem on quadratic differential forms that start our interest 
in combinatorial identities; in Section~\ref{identities}, we present some 
combinatorial formulae involving the binomial coefficients that, to our 
knowledge, are new; in Section~\ref{isotopy1} we apply these formulae to prove 
that some quadratic differential forms are isotopic; and in the last section  
we consider an alternative approach to prove isotopy. 

\section{Problem on Hessian Topology}\label{prelim}
Let us consider $H^n[x,y] \subset \mathbb R[x,y]$ the set of real homogeneous
polynomials of degree $n\geq 1\,$ in two variables. The graph of any $f \in H^n[x,y]$ 
contains the origin of $\mathbb R^3$. The polynomial $f$
is called {\it hyperbolic {\rm(}elliptic{\rm)}}
if its graph is a surface with only hyperbolic (elliptic) points off the origin.

The second fundamental form $II_{f}$ of a hyperbolic homogeneous polynomial $f$
defines two asymptotic lines at each point of $\mathbb R^{2^*}$.
Moreover, it defines two continuous asymptotic fields of lines without singularities 
on $\mathbb R^{2^*}$ that extend to the origin with a singularity. 
Since  both of these fields of lines are transverse at each point of $\mathbb R^{2^*}$,
their indexes at the origin coincide. 
Consequently, this index will be called 
{\it the index of the field of asymptotic lines at the origin},
and it will be denoted by $i_0(II_f)$.

A {\it hyperbolic isotopy} between two smooth hyperbolic quadratic 
differential forms $\omega$ and $\delta$
on $\mathbb R^{2^*}$ that extend themselves to the origin with a singularity is
a smooth map
$$\Psi: \mathbb R^{2^*} \times [0,1]  \rightarrow {\cal Q},\ \ \ (x,y,t) \mapsto \Psi_t(x,y),$$
\noindent where ${\cal Q}$ is the space of real quadratic forms on the plane and the following conditions hold:
$\Psi_0(x,y)= \omega(x,y)$, $\Psi_1(x,y)= \delta(x,y)$ and $\Psi_t(x,y)$ 
is a smooth hyperbolic quadratic differential form on 
$\mathbb R^{2^*}$, which extends at the origin with a singularity. In this case, 
we will say that {\it $\omega$
and $\delta$ are hyperbolic isotopic}.

If the second fundamental forms of two hyperbolic homogeneous polynomials of degree $n$
are hyperbolic isotopic, then  
the indexes of their fields of asymptotic lines at the origin coincide. 

The subset of $H^n[x,y]$, constituted by hyperbolic polynomials, 
is a topological subspace of $\mathbb R[x,y]$, denoted by $Hyp(n)$.
The topology of this space has been studied 
as part of the subject, known as Hessian Topology, 
introduced in \cite{arn1,panov,arn2}, and named 
by V. I. Arnold \cite{geoast}, \cite[problems 2000-1, 2000-2, 2001-1, 2002-1]{arn3}. In fact, Arnold stated the following conjecture \cite[p.1067]{geoast}, \cite[p.139]{arn3}:

\begin{description}
\item[ ] {\it ``The number of connected components of the space of hyperbolic \newline homogeneous
polynomials of degree $n$ increases as $n$ increases \quad \newline
 {\rm (}at least as a linear function of $n${\rm)}."}
\end{description}

In relation to this, he 
proved as an application of his characterization of hyperbolic polynomials in polar coordinates \cite[p.1031]{geoast},  that the homogeneous polynomials of degree $n$

$$ f(x,y) = (x^2+y^2)^{\frac{n-m}{2}} \mbox{Re} (x+i y)^m 
$$
are hyperbolic if $\, m\leq n < m^2 \,$ and $\,n-m\,$ is even. Moreover, he obtains that  
$\,i_0 (II_f) = \frac{2-m}{2}$.
This family of polynomials play a fundamental roll in finding the connected components 
mentioned in the conjecture. In the present article,
we prove as an application of our combinatorial formulae that 
the second fundamental form $II_f$ of $f$ is hyperbolic isotopy to the second fundamental form of $P(x,y):= \mbox{Re} (x+i y)^m$,  Theorem~\ref{propindices}.

\section{Combinatorial Identities}\label{identities}

 Let us start by presenting the combinatorial identities that we talk about. These identities occur naturally in the development of a direct proof of Theorem~\ref{propindices}, 
 but we decided to isolate
  them from their original context in order to be more general. 

\begin{theorem}\label{lemaA}
Let $m\geq 2\,$ be a natural number. For each integer number $\,0\leq j\leq \frac{m-1}{2}\,$ the following expression is fulfilled:
{\footnotesize
\begin{eqnarray}
(-1)^{j} \left[
\left(
\begin{array}{c}
m-1\\ 2j\\
\end{array}
\right) +
 \sum_{k=0}^{j-1} \left[
\left(
\begin{array}{c}
m-1\\ 2k\\
\end{array}
\right)
\left(
\begin{array}{c}
m-1\\ 2j-2k\\
\end{array}
\right) - \right.\right.\hskip 5.1cm \nonumber\\
\left.\left.\left(
\begin{array}{c}
m-1\\ 2k+1\\
\end{array}
\right)
\left(
\begin{array}{c}
m-1\\ 2j-2k-1\\
\end{array}
\right)
\right]\right] =
\left(
\begin{array}{c}
m-1\\ j\\
\end{array}
\right). \hskip 0.5 cm \label{coefA}
\end{eqnarray}}
\end{theorem}

\begin{theorem}\label{lemaBC}
If $m\geq 2\,$ is an even natural number, then for each integer number $\,0\leq j\leq \frac{m-2}{2}\,$ the following expressions are true:
{\footnotesize
\begin{eqnarray}
(-1)^{\frac{m}{2}+j-1}  \left[
- \left(
\begin{array}{c}
m-1\\ 2j-1\\
\end{array}
\right) +
 \sum_{k=0}^{\frac{m}{2}-j-1} \left[
\left(
\begin{array}{c}
m-1\\ 2k+2j\\
\end{array}
\right)
\left(
\begin{array}{c}
m-1\\ 2k+1\\
\end{array}
\right) - \right.\right.\hskip 3.36cm\nonumber\\
\left.\left.
\left(
\begin{array}{c}
m-1\\ 2k\\
\end{array}
\right)
\left(
\begin{array}{c}
m-1\\ 2k+2j-1\\
\end{array}
\right)
\right]\right]= \left(
\begin{array}{c}
m-1\\ j+\frac{m}{2}-1\\
\end{array}
\right).\hskip 0.5cm \label{coefC}\\
\nonumber \\
(-1)^{\frac{m}{2}}\left(
1-m + \sum_{k=0}^{\frac{m}{2}-2}
\left(
\begin{array}{c}
m-1\\ 2k+1\\
\end{array}
\right)
\left[
\left(
\begin{array}{c}
m-1\\ 2k+2\\
\end{array}
\right)-
\left(
\begin{array}{c}
m-1\\ 2k\\
\end{array}
\right)\right]\right)=
\left(
\begin{array}{c}
m-1\\ \frac{m}{2}\\
\end{array}
\right). \hskip 1.27cm \label{coefB}
\end{eqnarray}}
\end{theorem}

\noindent {\bf Proof of Theorem \ref{lemaA}}

Using the formula 
\begin{eqnarray}\label{fabsorcion}
(m-k) \left(\begin{array}{c}
m\\ k\\
\end{array}\right)= (k+1) \left(\begin{array}{c}
m\\ k+1\\
\end{array}\right),
\end{eqnarray}
several times, the expression (\ref{coefA}) is equivalent to the expression:
{\small
\begin{eqnarray}\label{Asimplificada}
(-1)^{j} \left[
\left(
\begin{array}{c}
m-1\\2j\\
\end{array}
\medskip
\right) +
\sum_{k=0}^{j-1}\frac{4k-2j+1}{m}
\left(
\begin{array}{c}
m\\2k+1\\
\end{array}
\right)
\left(
\begin{array}{c}
\medskip
m\\2j-2k\\
\end{array}
\right)\right]
=
\left(
\begin{array}{c}
m-1\\j\\
\end{array}
\right).
\end{eqnarray}}
By using the formula of the alternating sum of consecutive binomial coefficients,
\begin{eqnarray}\label{alternante}
(-1)^r \left(
\begin{array}{c}
\medskip
m-1\\ r\\
\end{array}
\right) = \sum_{k=0}^{r} (-1)^k
\left(
\begin{array}{c}
m\\ k\\
\end{array}
\right),
\end{eqnarray}
the expression (\ref{Asimplificada}) can be written as:
{\small
\begin{eqnarray*}
(-1)^{j} \sum_{k=0}^{j-1}\frac{4k-2j+1}{m}
\left(
\begin{array}{c}
m\\2k+1\\
\end{array}
\right)
\left(
\begin{array}{c}
\medskip
m\\2j-2k\\
\end{array}
\right)
=  \sum_{k=1}^{j} (-1)^{k+1}
\left(
\begin{array}{c}
m\\ j+k\\
\end{array}
\right).
\end{eqnarray*}}

Now, the sum on the right-hand side is divided into the sums of the even and odd terms. The sum on the left-hand side is divided into the sum of the first $\frac{j}{2}$ terms, if $j$ is even ( $\frac{j+1}{2}$ if $j$ is odd), and the remaining terms, but the terms are listed in reverse order. 

{\small
\begin{eqnarray*}
\sum_{r=1}^{\frac{j}{2}}\frac{1-4r}{m}
\left(
\begin{array}{c}
m\\j-2r+1\\
\end{array}
\right)
\left(
\begin{array}{c}
\medskip
m\\j+2r\\
\end{array}
\right)
+  \sum_{r=0}^{\frac{j}{2}-1} \frac{4r+1}{m}
\left(
\begin{array}{c}
m\\ j+2r+1\\
\end{array}
\right)
\left(
\begin{array}{c}
m\\ j-2r\\
\end{array}
\right) = \nonumber\\
= \sum_{r=1}^{\frac{j}{2}}(-1)
\left(
\begin{array}{c}
m\\j+2r\\
\end{array}
\right)
+  \sum_{r=0}^{\frac{j}{2}-1} 
\left(
\begin{array}{c}
m\\ j+2r+1\\
\end{array}
\right). \hskip 2.5cm
\end{eqnarray*}}

Then, by associating corresponding terms of both sides,
we obtain the equation

\begin{eqnarray}\label{formulafinal}
\sum_{k=1}^{j} \left[ (-1)^{k}
\left(
\begin{array}{c}
m\\j+k\\
\end{array}
\right)
\left( 1 + \frac{1-2k}{m}
\left(
\begin{array}{c}
m\\ j-k+1\\
\end{array}
\right)\right)\right] = 0.
\end{eqnarray}
We present now a proof of (\ref{formulafinal}) by means of recurrence
relations techniques. It was given by C. Merino-L\'opez \cite{cmerino}.
By using the alternating sum (\ref{alternante}) we have
{\scriptsize
\begin{eqnarray*}
m \sum_{k=0}^{j} (-1)^k
\left(
\begin{array}{c}
m\\j+k\\
\end{array}
\right)
= 
\left\{
\begin{array}{c}
(m-2j)\left(
\begin{array}{c}
m\\2j\\
\end{array}
\right)
- (m-j) \left(
\begin{array}{c}
m\\j\\
\end{array}
\right) \quad\mbox{{\normalsize if $j$ is even,}}\\
-(m-2j)\left(
\begin{array}{c}
m\\2j\\
\end{array}
\right)
- (m-j) \left(
\begin{array}{c}
m\\j\\
\end{array}
\right) \quad\mbox{{\normalsize if $j$ is odd.}}\\
\end{array}
\right.
\end{eqnarray*}}
By formula (\ref{fabsorcion}), the last equality becomes

{\scriptsize
\begin{eqnarray*}
m \sum_{k=0}^{j} (-1)^k
\left(
\begin{array}{c}
m\\j+k\\
\end{array}
\right)
= 
\left\{
\begin{array}{c}
(2j+1)\left(
\begin{array}{c}
m\\2j+1\\
\end{array}
\right)
- (j+1) \left(
\begin{array}{c}
m\\j+1\\
\end{array}
\right) \quad\mbox{{\normalsize if $j$ is even,}}\\
-(m-2j)\left(
\begin{array}{c}
m\\2j\\
\end{array}
\right)
- (m-j) \left(
\begin{array}{c}
m\\j\\
\end{array}
\right) \quad\mbox{{\normalsize if $j$ is odd.}}\\
\end{array}
\right.
\end{eqnarray*}}
Replacing the last equality in (\ref{formulafinal}) we obtain
\begin{eqnarray*}
\sum_{k=1}^{j+1}(-1)^{k+1} (2k-1)
\left(
\begin{array}{c}
m\\ k+j\\
\end{array}
\right)
\left(
\begin{array}{c}
m\\ j-k+1\\
\end{array}
\right)
=  (j+1) \left(
\begin{array}{c}
m\\ j+1\\
\end{array}
\right).
\end{eqnarray*}
Denote by $\,T(m,j)\,$ the left-hand side of the last expression. Using
 Stifel's identity,
{\scriptsize $\,\left(
\begin{array}{c}
m\\j\\
\end{array}
\right) = \left(
\begin{array}{c}
m-1\\j\\
\end{array}
\right) + \left(
\begin{array}{c}
m-1\\j-1\\
\end{array}
\right)$},
we verify that $\,T(m,j)\,$ satisfies
the recurrence relation
\begin{eqnarray}\label{recurrenciaT}
T(m,j) &=& T(m-1,j) + T(m-1,j-1) + \left(
\begin{array}{c}
m-1\\ j\\
\end{array}
\right)^2 + \nonumber\\
&&+ \,\, 2
\sum_{k=1}^{j}(-1)^{k}
\left(
\begin{array}{c}
m-1\\ j-k\\
\end{array}
\right)
\left(
\begin{array}{c}
m-1\\ j+k\\
\end{array}
\right).
\end{eqnarray}
Now, applying again  Stifel's formula to the function
\begin{eqnarray*}
F(m-1,j) = \left(
\begin{array}{c}
m-1\\ j\\
\end{array}
\right)^2 + 2
\sum_{k=1}^{j}(-1)^{k}
\left(
\begin{array}{c}
m-1\\ j-k\\
\end{array}
\right)
\left(
\begin{array}{c}
m-1\\ j+k\\
\end{array}
\right),
\end{eqnarray*}
it can be verified that $F$ satisfies the recurrence relation
\begin{eqnarray}\label{recurrenciaF}
F(m,j) = F(m-1,j) + F(m-1,j-1).
\end{eqnarray}
Stifel's identity establishes that the binomial coefficients also satisfy (\ref{recurrenciaF}).
Because  {\scriptsize $\,\left(
\begin{array}{c}
m\\j\\
\end{array}
\right)$}\, and $F(m,j)$ satisfy the same initial conditions,
we conclude that  $\, F(m,j) =$ {\scriptsize$\left(
\begin{array}{c}
m\\j\\
\end{array}
\right)$}.
So, the recurrence relation (\ref{recurrenciaT}) becomes
$$ T(m,j) - T(m-1,j) - T(m-1,j-1) =
\left(
\begin{array}{c}
m-1\\j\\
\end{array}
\right).$$
But, this relation is also fulfilled by
$\, (j+1) ${\scriptsize $\left(
\begin{array}{c}
m\\j+1\\
\end{array}
\right)$} with the same initial conditions, thus 
$ T(m,j) = (j+1) ${\scriptsize
$\left(\begin{array}{c}
m\\j+1\\
\end{array}\right)$} and~(\ref{formulafinal}) is proved.
\hfill$\Box$
\normalsize
\bigskip

\noindent {\bf Proof of Theorem \ref{lemaBC}}

Note that using formula (\ref{fabsorcion}) on the left side of $(\ref{coefC})$, it results in

{\small
\begin{eqnarray}\label{Csimplificado}
(-1)^{\frac{m}{2}+j-1} \left[
\sum_{k=0}^{\frac{m}{2}-j}\frac{m-4k-2j-1}{m}
\left(
\begin{array}{c}
m\\ 2k+1\\
\end{array}
\right)
\left(
\begin{array}{c}
m\\ 2k+2j\\
\end{array}
\right)\right]
=
\left(
\begin{array}{c}
m-1\\ j+\frac{m}{2}-1\\
\end{array}
\right).
\end{eqnarray}}

\medskip
\noindent Because $m$ is even, we make the substitutions $\,\, m = 2r\,$ and $\,\, j = r-n\,$ 
to obtain
{\small
\begin{eqnarray*}
(-1)^{2r-n-1} \left[
\sum_{k=0}^{n}\frac{(2n-4k-1)}{2r}
\left(
\begin{array}{c}
2r\\ 2k+1\\
\end{array}
\right)
\left(
\begin{array}{c}
2r\\ 2k+2r-2n\\
\end{array}
\right)\right]
=
\left(
\begin{array}{c}
2r-1\\ 2r-n-1\\
\end{array}
\right).
\end{eqnarray*}}

\medskip
\noindent
By using the {\it symmetry identity}
\,{\scriptsize $ \left(
\begin{array}{c}
a\\ b\\
\end{array}
\right)
= \left(
\begin{array}{c}
a\\ a-b\\
\end{array}
\right)\,$} 
(see \cite{graham}) in the last expression, and 
changing back $\, 2r=m\,$, the result is

\begin{eqnarray}\label{f3}
(-1)^{n} \left[
\sum_{k=0}^{n}\frac{(4k-2n+1)}{m}
\left(
\begin{array}{c}
m\\ 2k+1\\
\end{array}
\right)
\left(
\begin{array}{c}
m\\ 2n-2k\\
\end{array}
\right)\right]
=
\left(
\begin{array}{c}
m-1\\ n\\
\end{array}
\right).
\end{eqnarray}

\medskip
\noindent Note that the term $\,k=n\,$ on the left-hand side sum equals
{\scriptsize $\left(
\begin{array}{c}
m-1\\ 2n\\
\end{array}
\right)$}. 
Then, expression (\ref{f3}) is (\ref{Asimplificada}), and thus, equation (\ref{coefC})
is proved.
\bigskip

We finish with a proof of equation $(\ref{coefB})$.
The following expression is derived from formula $(\ref{fabsorcion})$:
\begin{eqnarray*}
\left( \begin{array}{c}
m-1\\ 2k+2\\
\end{array}
\right)-
\left(
\begin{array}{c}
m-1\\ 2k\\
\end{array}
\right) =
\frac{m-4k-3}{m}
\left(
\begin{array}{c}
m\\ 2k+1\\
\end{array}
\right)
\left(
\begin{array}{c}
m\\ 2k+2\\
\end{array}
\right).
\end{eqnarray*}
Then, equation $(\ref{coefB})$ is equivalent to
\begin{eqnarray}\label{f2}
(-1)^{\frac{m}{2}}\left(
\sum_{k=0}^{\frac{m}{2}-1}
\frac{(m-4k-3)}{m}
\left(
\begin{array}{c}
m\\ 2k+1\\
\end{array}
\right)
\left(
\begin{array}{c}
m\\ 2k+2\\
\end{array}
\right)
\right)=
\left(
\begin{array}{c}
m-1\\ \frac{m}{2}\\
\end{array}
\right).
\end{eqnarray}
But, when we replace $\, j=1$ in (\ref{Csimplificado}),
we retrieve (\ref{f2}).
\hfill$\Box$

\section{An application: isotopic  quadratic forms}\label{isotopy1}

In the following analysis, we consider a polynomial $f \in \mathbb R[x,y]$   
as a Hamiltonian function with
Hamiltonian vector field $\nabla f= (f_y, -f_x)$ 
on $\mathbb R^2$. The field of Hessian matrices
$Hess f =
\left(\begin{array}{cc}
f_{xx}& f_{xy}\\
f_{xy}& f_{yy}
\end{array}
\right)$ determines at each point $p \in \mathbb R^2$ a bilinear form. 
That is, 
$$Hess f_p:\mathbb R^2  \times \mathbb R^2 \rightarrow \mathbb R $$
$$Hess f_p(X,Y)= X (Hess f(p)) Y^t,$$
\noindent where $X, Y \in \mathbb R^2$ and the index $t$, means the transpose of the vector $Y$.

Thus, for any homogeneous polynomial $P \in H^{m}[x,y]$ we define
the following application:
\begin{eqnarray*}
\nabla P Hess P :& H^u[x,y] \rightarrow H^{2m+u-4}[x,y],\\  
& Q \mapsto \nabla P{HessP} \nabla Q^t.
\end{eqnarray*}

\noindent A straightforward computation implies that
\begin{eqnarray*}
\nabla P{HessP} \nabla Q^t= P_{xx}P_y Q_y + P_{yy} P_x Q_x -P_{xy}(P_x Q_y + P_y Q_x).
\end{eqnarray*}

\noindent The following inequality plays an important role in the proof 
of Theorem~\ref{propindices}. 
\begin{eqnarray} \label{uno}
\nabla P{HessP} \nabla Q^t(p) \leq 0,\ \ p \in \mathbb R^2.
\end{eqnarray}

We inform the reader that there are well-known techniques for dealing with the aforementioned inequality, and we present this alternative approach in the last section. But the purpose of this current paper is to highlight the combinatorial identities that we found, and their possible use in algebraic and differential geometry or elsewhere.

\begin{theorem}\label{propindices}
Let $P(x,y)=\mbox{Re} (x+i y)^m$ and $Q(x,y)=(x^2+y^2)^k$ 
be homogeneous polynomials of degree $m\geq 2$ with  $\, 2k\geq 2$ and $\, m\leq n <  m^2$. 
Then, the quadratic form $II_f$ is hyperbolic isotopic to the quadratic form $II_P$,
where $f=Q P$.
\end{theorem}

\noindent {\bf Proof.} 
A direct computation implies that the second fundamental form of the product $f=QP$ has the expression:
\begin{eqnarray}
\label{segundaformaproducto}
II_{PQ}= PII_Q + 2 dPdQ  + QII_P,
\end{eqnarray}
\noindent where $dP= P_x dx + P_y dy$, $dQ=Q_xdx + Q_y dy$ and $dPdQ$ is the quadratic differential form defined by the product of $dP$ and $dQ$.

Let us denote by $\omega$ the quadratic form $ 2 dP dQ(x,y) + Q(x,y)II_P(x,y)$.  We recall that {\it the discriminant of a quadratic form, $ a dx^2 + 2 b dx dy + c dy^2$,} is $\, b^2- ac.$ Firstly, we shall prove that $\omega$ is a hyperbolic quadratic form on $\mathbb R^{2^*}$, that is, that its discriminant $\Delta_\omega$ is positive on $\mathbb R^{2^*}$. After some computations, we have that
$$\Delta_\omega = -Q^2 \mbox{det}(Hess P) + \frac{1}{4}(P_x Q_y -P_y Q_x)^2 - 
2 Q(\nabla P{HessP} \nabla Q^t).$$
On one hand, the term $-Q^2 det(Hess P)$ is negative off the origin since $P$ is a
hyperbolic polynomial and $Q$ is positive off the origin.
On the other hand, the term $\Delta_{dPdQ}= \frac{1}{4}(P_x Q_y - P_y Q_x)^2$ is nonnegative on the plane. Finally, the proof follows from  Lemma \ref{condPQpar}.

Now, we shall prove that $\omega$ and
$\omega + \delta = II_f$ are hyperbolic isotopic, where $\delta$ is the quadratic form $P(x,y)II_Q(x,y)$. Consider the isotopy
$$ \Phi_t(x,y) = \omega(x,y) + t \delta(x,y). $$

Because $\,\omega(x,y)=\omega_1 dx^2 + 2 \omega_2 dx dy + \omega_3 dy^2$
and $\delta(x,y)= \delta_1 dx^2 + 2 \delta_2 dx dy + \delta_3 dy^2$ we have that 
the discriminant of
$ \Phi_t(x,y)$ is
$$   \Delta \Phi_t = \omega_2^2 -  \omega_1 \omega_3 + t (2 \omega_2\delta_2- \omega_1\delta_3 \omega_3\delta_1) 
+ t^2 (\delta_2^2 - \delta_1 \delta_3).$$
Note that
$\omega_2^2 -  \omega_1 \omega_3 +  (2 \omega_2\delta_2- \omega_1\delta_3 \omega_3\delta_1) + (\delta_2^2 -  \delta_1 \delta_3)>0\,$ and $\,\omega_2^2 -  \omega_1 \omega_3 > 0\ $ 
 on $\mathbb R^{2^*}$ since $\omega$  and  $\omega +\delta$ are hyperbolic 
on $\mathbb R^{2^*}$. So, for $t \in [0,1]$ we have that the discriminant $ \Delta \Phi_t$
 is positive on $\mathbb R^{2^*}$.

The next step is to prove that the quadratic differential forms 
$Q II_p  + 2  dP dQ$ and $Q II_p$
are hyperbolic isotopic. To do that, consider the isotopy 
$\Psi_t(x,y)= Q II_p + 2t dP dQ (x,y)$, where  $\,t \in [0,1]$.
We can see  that the discriminant of the quadratic differential form 
$Q II_p  + 2t  dP dQ (x,y)$ is
$$\Delta_{\Psi_t}= -Q^2 det(Hess P) + t^2 
\Delta_{dPdQ} - 2t Q(\nabla P{HessP} \nabla Q^t),$$ 
which is positive on $\mathbb R^{2^*}$ by the next Lemma.
\hfill$\Box$

\begin{lemma} \label{condPQpar}
Let  $m, k \in \mathbb Z$ such that $m \geq 2$ and  $k \geq 1$.
Then, the homogeneous polynomials  $P, Q$ of degree $m$ and $2k$ as previously defined 
satisfy inequality $(\ref{uno})$.
\end{lemma}

\noindent {\bf Proof of Lemma \ref{condPQpar}.}
We only present the case when $m$ is even. The odd case is analogous.
In order to prove that inequality $(\ref{uno})$ holds for the 
polynomials $P$ and $Q$ we consider the polynomial expression 
$\nabla P {HessP} (\nabla Q)^t$ and prove that
$$ P_x(Q_y P_{xy}  - Q_x P_{yy}) + P_y (Q_x P_{xy} - Q_y P_{xx}) = 2\,k \, m^2 (m-1) (x^2+y^2)^{k+m-2}.$$
\noindent Since $m$ is even, a straightforward computation shows that
$$Q_y P_{xy} - Q_x P_{yy} =2k (x^2+y^2)^{k-1}
\left[
\sum_{j=0}^{\frac{m}{2}-1} (-1)^j
\frac{ (m-1) m!}{(2j)! (m-2j-1)!}
x^{m-2j-1} y^{2j} \right].$$

\noindent Now, we multiply both sides of the last expression by $P_x$. The
product $\,P_x \,(Q_y P_{xy} - Q_x P_{yy})\,$ equals
\medskip
\begin{eqnarray*}
 2k (x^2+y^2)^{k-1} m^2 (m-1) \left[
\sum_{j=0}^{\frac{m}{2}-1} (-1)^j
\left(
\begin{array}{c}
m-1\\ 2j\\
\end{array}
\right)
x^{m-2j-1} y^{2j} \right]^2.
\end{eqnarray*}

\noindent Developing the squared factor of the last expression we have
$$P_x \,(Q_y P_{xy} - Q_x P_{yy}) \,=\, 2\,k \,(x^2+y^2)^{k-1} m^2 (m-1) \hskip 4.3cm$$
$$ \quad\left[
x^{2m-2} \,+ \sum_{j=1}^{\frac{m}{2}-1} \left(\sum_{k=0}^j (-1)^j
\left(
\begin{array}{c}
m-1\\ 2k\\
\end{array}
\right)
\left(
\begin{array}{c}
m-1\\ 2j-2k\\
\end{array}
\right)\right)
x^{2m-2j-2} y^{2j} \,\, +  \right.$$
{\small
\begin{eqnarray}\label{ecua1}
\left.\qquad\sum_{j=1}^{\frac{m}{2}-1} 
\left(\sum_{k=0}^{\frac{m}{2}-j-1} (-1)^{\frac{m}{2}+j-1}
\left(
\begin{array}{c}
m-1\\ 2k+2j\\
\end{array}
\right)
\left(
\begin{array}{c}
m-1\\ m-2k-2\\
\end{array}
\right)\right)
x^{m-2j} y^{m+2j-2}\right].
\end{eqnarray}}

\noindent By doing a similar computation for $ P_y (Q_x P_{xy} - Q_y P_{xx})$ we obtain that 

$$P_y \,(Q_x P_{xy} - Q_y P_{xx}) \,=\, 2\,k \,(x^2+y^2)^{k-1} m^2 (m-1)\hskip 4.2cm$$
{\small
\begin{eqnarray*}
\qquad\left[ y^{2m-2} + \sum_{j=1}^{\frac{m}{2}} \left(\sum_{k=0}^{j-1} (-1)^{j+1}
\left(
\begin{array}{c}
m-1\\ 2k+1\\
\end{array}
\right)
\left(
\begin{array}{c}
m-1\\ 2j-2k-1\\
\end{array}
\right)\right)
x^{2m-2j-2} y^{2j} \,\,+ \right.
\end{eqnarray*}}
{\small
\begin{eqnarray}\label{ecua2}
\left.\qquad\sum_{j=2}^{\frac{m}{2}-1} \left(\sum_{k=0}^{\frac{m}{2}-j} (-1)^{\frac{m}{2}+j}
\left(
\begin{array}{c}
m-1\\ 2k+2j-1\\
\end{array}
\right)
\left(
\begin{array}{c}
m-1\\ m-2k-1\\
\end{array}
\right)\right)
x^{m-2j} y^{m+2j-2}\right].
\end{eqnarray}}

\noindent By adding the expressions  (\ref{ecua1}) and (\ref{ecua2}) we obtain
 \begin{eqnarray}\label{suma-simpli} \nonumber
P_x \,(Q_y P_{xy} - Q_x P_{yy}) + P_y \,(Q_x P_{xy} - Q_y P_{xx}) = \hskip 4.5cm\\ \nonumber
= 2\,k \,(x^2+y^2)^{k-1} m^2 (m-1)\Bigg[ x^{2m-2}\quad\hskip 2.4cm\\ \nonumber
 + \bigg(\sum_{j=1}^{\frac{m}{2}-1} A(j) x^{2m-2j-2} y^{2j}\bigg) + B x^{m-2} y^{m}
 \hskip 2.4cm \\ 
+ \left. \bigg(\sum_{j=2}^{\frac{m}{2}-1} C(j) x^{m-2j} y^{m+2j-2}\bigg) + y^{2m-2} \right].
\hskip 2.3cm 
\end{eqnarray}

\noindent By replacing (\ref{coefA}), (\ref{coefB}) and (\ref{coefC}) in (\ref{suma-simpli}) we conclude that
\bigskip

\noindent $P_x \,(Q_y P_{xy} - Q_x P_{yy}) + P_y \,(Q_x P_{xy} - Q_y P_{xx}) =$
{\small
\begin{eqnarray*}
= 2\,k \,(x^2+y^2)^{k-1} m^2 (m-1)
\left[\sum_{j=0}^{\frac{m}{2}-1} \left(
\begin{array}{c}
m-1\\ j\\
\end{array}
\right) x^{2m-2j-2} y^{2j}  \,\,+ \right.\hskip 0.9cm\\
\left. \hskip 2.3cm +
\left(
\begin{array}{c}
m-1\\ \frac{m}{2}\\
\end{array}
\right)
x^{m-2} y^{m} + \sum_{r=\frac{m}{2}+1}^{m-2}
\left(
\begin{array}{c}
m-1\\ r\\
\end{array}
\right)
x^{2m-2r-2} y^{2r}+ \,y^{2m-2}\right].
\end{eqnarray*}}

\noindent By collecting all the terms  inside  the square brackets we have
\bigskip

\noindent $P_x \,(Q_y P_{xy} - Q_x P_{yy}) + P_y \,(Q_x P_{xy} - Q_y P_{xx}) =$
{\small
\begin{eqnarray*}
\hskip 3.1cm = 2\,k \,(x^2+y^2)^{k-1} m^2 (m-1)\left[
\sum_{j=0}^{m-1}
\left(
\begin{array}{c}
m-1\\ j\\
\end{array}
\right)
x^{2m-2j-2} y^{2j}\right].
\end{eqnarray*}}

\noindent  Note that the expression  inside the square brackets equals
the binomial $\,(x^2+y^2)^{m-1}$.  So, we conclude
$$P_x(Q_y P_{xy} - Q_x P_{yy})+ P_y(Q_x P_{xy} - Q_y P_{xx}) =
2\,k \, m^2 (m-1)
\left(x^2+y^2\right)^{k+m-2}
\Box$$

\section{Alternative Approach}\label{isotopy2}
Let us present a much faster way of obtaining inequality $(\ref{uno})$.
Consider  Euler's Lemma:

Let $f\in H^n[x,y]$ be a homogeneous polynomial of degree $n$. Then,
$$ n f(x,y)= x f_x(x,y) + y f_y(x,y).$$ 
By doing a recursive use of this Lemma we obtain that
$$ P_x(Q_y P_{xy}  - Q_x P_{yy}) + P_y (Q_x P_{xy} - Q_y P_{xx}) = \frac{n}{m-1} 
Q (Hess P).$$
So, $\nabla P {HessP} (\nabla Q)^t \leq 0$ since $P$ is a hyperbolic polynomial and
$Q$ is positive on $\mathbb R^{2^*}$.

However, as mentioned before, we wanted to used the combinatorial identities 
of Theorem~\ref{lemaA} and~\ref{lemaBC} as a way of putting combinatorics to the 
service of algebraic and differential geometry and with the hope that they can 
be used elsewhere.  
\bigskip

\noindent {\bf Acknowledgments}
We would like to thank C. Merino-L\'opez for his proof of
equation (\ref{formulafinal}). We also thank Francisco
Larri\'on for nice conversations on this subject.
Work partially supported by DGAPA-UNAM grant PAPIIT-IN111415, PAPIIT-N102413 and
PAPIIT-IN110803.

\noindent Adriana Ortiz-Rodr\'{\i}guez \newline 
Instituto de Matem\'aticas, Universidad Nacional Aut\'onoma de M\'exico \newline
Area de la Inv. Cient., Circuito exterior C.U., M\'exico D.F 04510, M\'exico.\newline
e-mail: aortiz@matem.unam.mx\\

\noindent Federico S\'anchez-Bringas \newline
Depto. de Matem\'aticas, Fac. de Ciencias, Universidad Nacional Aut\'onoma de M\'exico,
Ciudad Universitaria, M\'exico D.F 04510, M\'exico.\newline
e-mail: sanchez@unam.mx

\end{document}